\documentclass{amsart}
\title [Tropical covers of curves]{Tropical covers of curves and their moduli spaces}
\author {Arne Buchholz and Hannah Markwig}
\address {Arne Buchholz, Universit\"at des Saarlandes, Fachrichtung
  Mathematik, Postfach 151150, 66041 Saarbr\"ucken, Germany}
\email {buchholz@math.uni-sb.de}
\address {Hannah Markwig, Universit\"at des Saarlandes, Fachrichtung
  Mathematik, Postfach 151150, 66041 Saarbr\"ucken, Germany}
\email {hannah@math.uni-sb.de}
\thanks {\emph {2010 Mathematics Subject Classification:} Primary 14T05, Secondary 14N35, 51M20}
\keywords {Tropical geometry, Hurwitz numbers, covers of curves}

\usepackage{amsfonts}
\usepackage{amssymb}
\usepackage{amscd}
\usepackage{amsmath}
\usepackage{graphicx}
\usepackage{color}
\allowdisplaybreaks[3]

\newcommand{\bbR}{\mathbb{R}}

\newcommand{\bbP}{\mathbb{P}}

\def\bary{\begin{array}} 
\def\eary{\end{array}} 
\def\ben{\begin{enumerate}} 
\def\een{\end{enumerate}}
\def\bit{\begin{itemize}} 
\def\eit{\end{itemize}}

\newcommand {\dunion}{\,\mbox {\raisebox{0.25ex}{$\cdot$} \kern-1.83ex $\cup$}
  \,}

\newcommand{\cC}{\mathcal{C}}
\newcommand{\CL}{\mathcal{L}}
\newcommand {\PP}{{\mathbb P}}
\newcommand {\RR}{{\mathbb R}}
\newcommand {\ZZ}{{\mathbb Z}}
\newcommand {\Perm}{{\mathbb S}}
\newcommand {\N}{{\mathbb N}}

\DeclareMathOperator{\trop}{trop}

\DeclareMathOperator{\Mg}{M_g^{trop}}

\DeclareMathOperator{\Mone}{M_1^{trop}}
\DeclareMathOperator{\val}{val}
\DeclareMathOperator{\Aut}{Aut}
\DeclareMathOperator{\br}{br}
\DeclareMathOperator{\mult}{mult}

\def\beq{\begin{equation}}                     %  
\def\eeq{\end{equation}}                       % 
\def\bea{\begin{eqnarray}}                     %         % 
\def\eea{\end{eqnarray}}

\theoremstyle{plain}

\newtheorem{thm}{Theorem}[section]
\newtheorem{lemma}[thm]{Lemma}
\newtheorem{prop}[thm]{Proposition}

\newtheorem*{conj*}{Conjecture}
\newtheorem{cor}[thm]{Corollary}
\newtheorem*{cor*}{Corollary}

\theoremstyle{definition}

\newtheorem{remark}[thm]{Remark}
\newtheorem{definition}[thm]{Definition}
\newtheorem{example}[thm]{Example}

\newcommand{\GIT}[1]{/\!\!/_{\kern-.2em #1 \kern0.1em}}

\begin {document}

\begin {abstract}
We define the tropical moduli space of covers of a tropical line in the plane as weighted abstract polyhedral complex, and the tropical branch map recording the images of the simple ramifications. Our main result is the invariance of the degree of the branch map, which enables us to give a tropical intersection-theoretic definition of tropical triple Hurwitz numbers. We show that our intersection-theoretic definition coincides with the one given in \cite{BBM10} where a Correspondence Theorem for Hurwitz numbers is proved. Thus we provide a tropical intersection-theoretic justification for the multiplicities with which a tropical cover has to be counted. Our method of proof is to establish a local duality between our tropical moduli spaces and certain moduli spaces of relative stable maps to $\PP^1$.
\end {abstract}

\maketitle

%%%%%%%%%%%%%%%%%%%%%%%%%%%%%%%%%%%%%%%%%%%%%%%%%%%%%%%%%%%%%%%%%%%%%%%%%%%%%%%

\section{Introduction}
Tropical geometry studies the geometry over the tropical semiring. It can be
viewed as a piece-wise linear degeneration of algebraic geometry which preserves
many properties but can be studied by combinatorial methods. Tropical geometry
has been particularly succesful for the study of enumerative geometry. Beginning
with Mikhalkin's Correspondence Theorem for the numbers of degree $d$ genus $g$
plane curves through $3d+g-1$ points in general position \cite{Mi03}, many
situations have been studied where an (algebraic) enumerative number agrees with
the corresponding tropical number. Such correspondence theorems make it possible
to study properties of the algebraic numbers (e.g.\ relations between them) by
means of tropical geometry. 
 
In enumerative geometry, a strategy to count is to find a suitable moduli space
parametrizing the objects to count, and then to determine the numbers as
intersection numbers on this moduli space. Likewise, in tropical geometry a
theory of moduli spaces and intersection theory has developed and many tropical
enumerative numbers can now be expressed as intersection numbers on an
appropriate moduli space (see e.g.\ \cite{CJM10}, \cite{GKM07}, \cite{Rau08}). 
\\
\\
Hurwitz numbers count covers of the projective line with fixed ramification
profile over given points. If we fix only two special ramification profiles and
simple ramification otherwise, we speak of double Hurwitz numbers. By matching a
cover with a monodromy representation, such a count is equivalent to choices of
$n$-tuples of elements of the symmetric group $S_d$ multiplying to the identity
element and acting transitively on the set $\{1 \ldots d\}$. The study of
Hurwitz numbers has provided a rich interplay between the combinatorics and
representation theory of the symmetric group and the geometry of covers of a
line.

Intersection-theoretically, Hurwitz numbers can be described as the degree of a
suitable branch map recording the images of the ramification points.
\\
\\
Tropical double Hurwitz numbers have been introduced in \cite{CJM10} by means of tropical intersection theory, as the degree of a tropical branch map recording the images of the simple ramifications, in analogy to the algebraic situation. Also, a correspondence theorem is proved stating the equality with the algebraic counterparts. Tropical double Hurwitz numbers have been useful in the study of the piece-wise polynomial structure of double Hurwitz numbers \cite{cjm:wcfdhn}. By analyzing the intersection-theoretic definition, one can see that tropical double Hurwitz numbers count covers of the simplest model of tropical $\bbP^1$, $\bbR\cup\{\pm \infty\}$, where the two special ramifications are imposed in terms of weight conditions of the ends mapping to $\pm \infty$ and the simple ramifications correspond to interior trivalent vertices of the source curve.
The moduli space of such tropical curves as a set consists of all such covers where the simple ramification points, i.e.\ the images of the trivalent vertices, are not fixed but can move around. When computing the degree of the branch map, we fix the images of all trivalent vertices.

In \cite{BBM10}, the tropical definition is generalized to arbitrary Hurwitz numbers and a correspondence theorem is proved.
Such more general tropical Hurwitz numbers play an important role e.g.\ in the tropical study of Zeuthen numbers \cite{BBM11}.
 However, the definition requires all the ramification to be at the ends of the target tropical curve, no simple ramification is allowed in the interior of the curve. With this restriction it does not make sense to consider tropical moduli spaces of such covers as they are just zero-dimensional --- there is no simple ramification to move. Thus \cite{BBM10} does not consider the tropical intersection-theoretic approach to Hurwitz numbers.
Consequently, the multiplicity with which a tropical cover has to be counted in the definition of \cite{BBM10} is modeled exactly right to satisfy the correspondence theorem, but it lacks an intrinsic justification within tropical geometry.

We close this gap by generalizing the definition of tropical covers of \cite{BBM10} allowing simple ramification also in the interior at trivalent vertices. We study the moduli space of such covers as an abstract weighted polyhedral complex and the branch map recording the images of the trivalent vertices. Our main theorem, Theorem \ref{thm}, states that the degree of this branch map is constant which enables us to redefine the tropical Hurwitz number as this degree, independently of the images that we fix for the trivalent vertices.

Our main theorem is formulated for the case of covers of a tropical line $\CL$ in the plane. This can be viewed as the building block for the general case of covers of a trivalent curve. To keep notations simple, we restrict the formulation to the case of covers of $\CL$.

Our result sheds light on the definition of multiplicity of a cover of \cite{BBM10}: we can analyze which factors are contributions coming from the weights of cells of the moduli space and which factors arise as local multiplicities of the branch map.

An interesting feature of our main result is that we use methods from algebraic geometry to prove it: using a dual graph construction, we match top-dimensional cells around a fixed cell of codimension one of our tropical moduli space with points in a one-dimensional moduli space of algebraic covers determined by the cell of codimension one. The degree of the tropical branch map can then be related to the degree of the branch map from the algebraic one-dimensional moduli space to $\bbP^1$, and the fact that it is constant follows from the fact that pull-backs of different points of $\bbP^1$ are equivalent.
Hence the methods of our paper complement the more common situation where tropical methods are used to derive results in algebraic geometry. 

For the case of genus zero, moduli spaces of covers similar to ours have been studied recently in \cite{GMO}. There, these spaces are embedded into a real vector space as a tropical variety (i.e.\ a balanced polyhedral complex). In the situations that are considered both in \cite{GMO} and in our paper, we show that the definitions of weights for the top-dimensional cells of the moduli spaces agree. 

For higher genus, moduli spaces of tropical curves are still an object of intense study \cite{CMV12}, \cite{ACP12}. The moduli spaces we consider can be mapped with the forgetful map forgetting the map and the target to the moduli space of tropical curves. The image of our moduli space under the forgetful map equals the set of curves that allow a cover of $\bbP^1$ satisfying given properties and can be thus viewed as a tropical analogue of a Hurwitz scheme. We are interested in relating these tropical Hurwitz schemes to their algebraic counterparts. 
\\
\\
Our paper is organized at follows. In section 2 we recall necessary definitions and construct our tropical moduli space of covers as an abstract weighted polyhedral complex. We compare our weights in the rational case to weights that recently appeared in a different construction of tropical moduli spaces. In section 3 we define the branch map, state our main result and deduce by comparing to the definition of \cite{BBM10} that our intersection-theoretic definition of tropical Hurwitz numbers yields the analogous algebraic Hurwitz numbers.
We then turn towards the proof of our main theorem. We first consider the main building blocks for our proof, namely one-dimensional moduli spaces and the duality of these spaces to algebraic moduli spaces of relative stable maps. Finally, we collect all partial results to finish the proof.

\subsection{Acknowledgements}
We thank the Deutsche Forschungsgemeinschaft for funding by DFG-grant MA 4797/1-1.
We thank Erwan Brugall\'{e}, Renzo Cavalieri, Andreas Gathmann and Dennis Ochse for helpful discussions.

\section{The moduli space of tropical covers of a line}
 We quickly recall basic definitions and fix notations. For more details on tropical curves and their morphisms, see e.g. \cite{BBM10, Cap12a, CJM10, Cha10}.

An \emph{(abstract, marked) tropical curve (with labelled vertices)} is a connected metric graph $\Gamma$ satisfying the following properties. A vertex is called a \emph{leaf} if it is one-valent and (inner) vertex otherwise. An edge $e$ is called \emph{end} and has length $l(e)=\infty$ if it is adjacent to a leaf, otherwise it is called a \emph{bounded edge} and has a length $l(e)\in \RR$. Each end is \emph{marked} by a number.
Each inner vertex $V$ is equipped with a number $g_V\in \N$ that we call the \emph{genus of the vertex}. If $\val(V)=2$ then $g_V\geq 1$. In addition, each inner vertex is \emph{labelled} with finitely many distinct numbers such that the disjoint union of all labels equals $\{1,\ldots,r\}$ for some $r$.

Two tropical curves are isomorphic (and will be identified in the following) if there exists an isomorphism of the underlying metric graphs preserving the genus, the labeling of the vertices, and the marking of the ends.

The number $g=b^1(\Gamma)+\sum_V g_V$ is called the \emph{genus} of the tropical curve $\Gamma$.

The \emph{combinatorial type} of a tropical curve is obtained by omitting the length data.

\begin{example}\label{exmp:tropCurve}
  Figure \ref{fig-curve} shows a genus $5$ tropical curve. The red numbers denote
the
genus on vertices, the black numbers are the edge lengths, vertex labels as well as end markings are left
out. All vertices without a red number have genus zero.

\begin{figure}
\begin{center}
 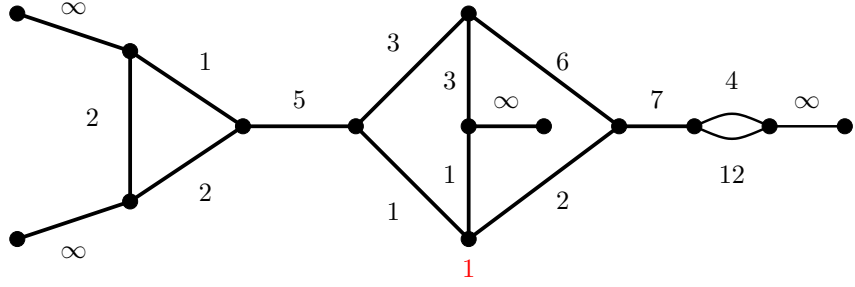
\end{center}
\caption{A tropical curve.}
\label{fig-curve}
\end{figure}
\end{example}

\begin{definition}[cf.\ \cite{BBM10}] A continous map $h:\Gamma\rightarrow \tilde\Gamma$ of tropical curves is called a \emph{morphism of tropical curves}, if
\begin{itemize}
 \item The image $h(e)$ for every edge $e$ of $\Gamma$ is contained in an edge of $\tilde\Gamma$.
 \item $h$ is integral affine-linear on each edge $e$, i.e.\ if we understand
$e$ as open interval $(0,l(e))$, then
$h_{\mid e}$ maps $t\in (0,l(e))$ to $w_et+a$ for some starting point $a$ on the image edge and some nonzero integer $w_e$ which is defined up to sign and called the \emph{weight} of $e$. 
 \item $h$ fulfills the \emph{balancing condition}:
For every vertex $V$ of $\Gamma$, let $\tilde{e}_1,\ldots,\tilde{e}_l$ be the edges
adjacent to $h(V)$ (if $h(V)$ is not a vertex of $\tilde\Gamma$ then consider it momentarily
as a $2$-valent vertex of genus $\tilde g_{h(V)}=0$ by subdividing the edge it is mapped to). For $i\in\{1,\ldots,l\}$
denote by $e^{(1)}_i,\ldots,e^{(k_i)}_i$ the edges adjacent to $V$ that
are mapped to $\tilde{e}_i$. Then
\begin{equation}
 \sum_{m=1}^{k_i}w_{e^{(m)}_i}=\sum_{m=1}^{k_j}w_{e^{(m)}_j} \label{balCond}
\end{equation}
for all $i,j$. The number in (\ref{balCond}) is called the \emph{local degree}
of $h$ at $V$ and is denoted by $\deg_hV$. For a point $a$ on an edge $e$ of $\Gamma$, we define the local degree to be equal to $\deg_ha=|w_e|$.
\item For every inner vertex $V$ of $\Gamma$, the \emph{Riemann-Hurwitz
condition}
 is fulfilled, i.e. 
\begin {equation}\label{eqn:RHnumber}
r_V:=(\val(V)+2g_V-2)-\deg_hV\cdot (\val(h(V))+2\tilde{g}_{h(V)}-2)\geq0.
\end {equation}
We will call $r_V$ the \emph{RH-number} of $V$. If $r_V$ is positive we
say that $V$ is \emph{ramified} or a \emph{ramification point}. In this case its image $h(V)$ is called \emph{branch point}.
If $r_V=1$ the vertex $V$ is a \emph{simple} ramification.
\item Each inner vertex $V$ has exactly $r_V$ different labels.
 \end{itemize} 
A \emph{(tropical) cover} of $\tilde \Gamma$ is a curve $\Gamma$ with a morphism $h:\Gamma \rightarrow \tilde \Gamma$ as above.
\end{definition}

Sometimes one also allows edges of weight $w_e=0$, i.e.\ edges which are contracted to a point. 
As contracted edges do not play a role for counting covers, we neglect them
here.

\begin{definition}
 Let $h:\Gamma \rightarrow \tilde \Gamma$ be a cover. The balancing
condition implies that for every point $\tilde V$ in the image the sum
\begin{equation}
 \sum_{V|h(V)= \tilde V} \deg_hV
\end{equation}
is the same. This number is called the \emph{degree} $\deg(h)$ of $h$.
\end{definition}

\begin{example}
 
Figure \ref{fig-cover} shows some local parts of tropical covers, resp.\ a map
which is not a tropical cover since it does not satisfy the RH-condition at the
vertex $V$. As usual, we do not show edge lengths, end markings and vertex labelings in the picture. The
red numbers denote the genus of the vertices, the blue numbers the weights of the
edges.

\begin{figure}
\begin{center}
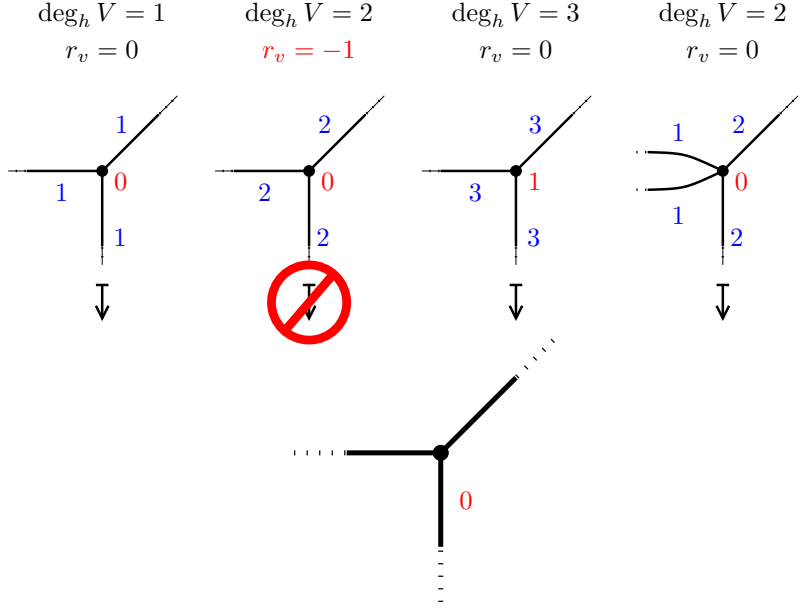\label{exmp:localCover}
\end{center}
\caption{Local pictures of tropical covers resp.\ a map which does not satisfy the RH-condition.}
\label{fig-cover}
\end{figure}
\end{example}

\begin{definition}\label{defn:automorphism}
 Let $h:\Gamma \rightarrow \tilde \Gamma$ be a cover. An \emph{automorphism} of
a cover $h:\Gamma \rightarrow \tilde \Gamma$ is an isomorphism $\phi:\Gamma\rightarrow\Gamma$, s.t. $\phi\circ h=h$. We
denote by $\Aut(h)$ the group of automorphisms of $h:\Gamma \rightarrow \tilde \Gamma$.
\end{definition}

We will focus our attention on covers of a tropical line in the tropical projective plane. To fix notation, we define:
\begin{definition}
Consider the abstract curve $\CL$ of a tropical line in the tropical projective plane, i.e.\ a curve with one vertex that we denote by $c$ and three ends adjacent to $c$ that we call $u$, $v$ and $w$ respectively (see Figure \ref{fig-L}).

\begin{figure}
 \begin{center}
 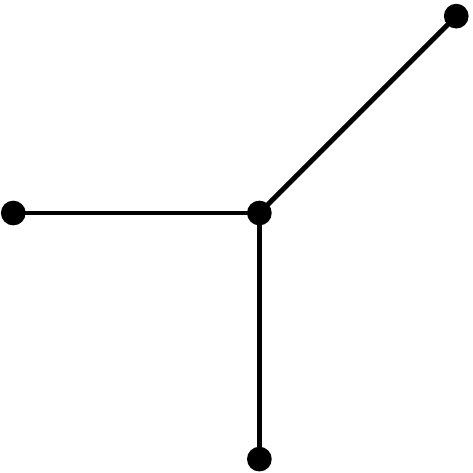
\end{center}\caption{The tropical line $\CL$.}\label{fig-L}
\end{figure}
\end{definition}

\begin{example}\label{exmp:cover}
 Figure \ref{fig-cover2} shows an example of a cover of $\CL$ by the curve $\cC$
of example \ref{exmp:tropCurve}.
We only mark the edge weights in blue, the other values can be deduced from \ref{exmp:tropCurve}. Note that the RH-numbers of the vertices mapping to $c$ are zero.
When drawing pictures of a cover of $\CL$ in the following, we will leave out the target $\CL$. The way we organize the picture of the curve mapping to $\CL$ indicates clearly which parts are mapped to $u$, $v$, $w$ and $c$, respectively.
\begin{figure}
\begin{center}
 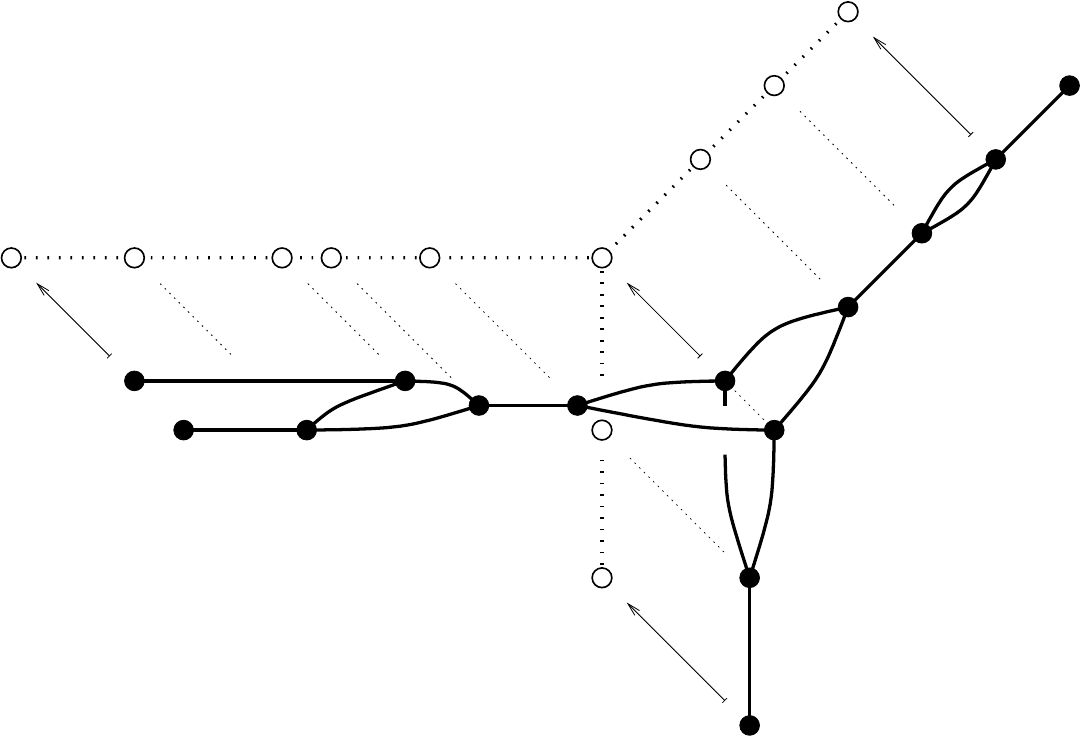
\end{center}
\caption{A tropical cover of $\CL$.}\label{fig-cover2}
\end{figure}

\end{example}

\begin{definition}
 The \emph{combinatorial type} of a cover $h:\Gamma \rightarrow \CL$ consists of the combinatorial type of $\Gamma$ together with the edge weights for all edges, and the set of vertices which is mapped to $c$.
\end{definition}

For a fixed combinatorial type $\alpha$, the set of covers of type $\alpha$ forms an open rational polyhedral cone that we call $D_\alpha$: we can vary the lengths of the bounded edges, but we cannot vary them independently since we need to cover $\CL$. Inside the open positive orthant of the real vector space of dimension the number of bounded edges, the conditions to have a cover of $\CL$ can be expressed as integral linear equations. We will see an example in \ref{exmp:equations}.

\begin{definition}\label{def-contraction}
Points on the boundary of the cone $\overline{D}_\alpha$ correspond to covers where some lengths of edges are shrunk to zero. We remove edges of zero lengths, identify their adjacent vertices and adjust the genus at vertices and their labels as follows:
Denote by $\Gamma'$ a connected subgraph of edges whose lengths go to zero. 
Let $V_1,\ldots,V_k$ be the vertices of $\Gamma'$.
Replace $\Gamma'$ by a vertex $V$ of genus $g_V=\sum_{i=1}^k g_{V_i}+b^1(\Gamma')$ and with the union of all labels of the $V_i$ as labels.

We call the new cover (resp.\ the new combinatorial type) obtained in this way a \emph{contraction} of $\alpha$ (resp.\ of a cover of type $\alpha$). 
\end{definition}
\begin{lemma}
 Definition \ref{def-contraction} is well-defined, i.e.\ a contraction
corresponding to a point on the boundary of the cone $\overline{D}_\alpha$ of a
combinatorial type $\alpha$  is indeed a cover.
\end{lemma}
\begin{proof}
Let $V$ be a new vertex replacing the connected subgraph $\Gamma'$ of edges going to zero. Assume $\Gamma'$ has $k$ vertices $V_1,\ldots,V_k$ and $E$ edges. Then $V$ has $\sum_i r_{V_i}$ labels and we have to see that it has RH-number $r_V=\sum_i r_{V_i}\geq 0$.
 Assume that for $l<k$ the vertices $V_1,\ldots,V_l$ are mapped to the center $c$ of $\CL$ while the $V_i$ with $i>l$ are mapped to a ray. Then the RH-numbers of the $V_i$ equal $r_{V_i}= \val(V_i)+2g_{V_i}-2-d_i$ if $i\leq l$ where $d_i$ denotes the local degree at $V_i$, and $r_{V_i}= \val(V_i)+2g_{V_i}-2$ else.
If $l=0$ all the $V_i$ as well as the new vertex $V$ must be mapped to the same ray and we have 
$r_V=\val(V)+2g_V-2$. If $l>0$, the new vertex $V$ must be mapped to $c$ and we have $r_V=\val(V)+2g_V-2-d$, where $d$ denotes the local degree at $V$. Obviously $d=\sum d_i$ in this case.
In any case we have \begin{align*}\sum_i r_{V_i} & =  \sum_{i=1}^k ( \val(V_i)+2g_{V_i}-2) -\sum_{i=1}^l d_i  \\ &=
\sum_{i=1}^k  \val(V_i) + 2 \sum_{i=1}^k g_{V_i} - 2k - \sum_{i=1}^l d_i \\ &=
\sum_{i=1}^k  \val(V_i) - 2 E + 2 \sum_{i=1}^k g_{V_i} +2(E-k+1) -2 - \sum_{i=1}^l d_i \\ & =
\val(V) + 2g_V -2 - \sum_{i=1}^l d_i = r_V,\end{align*}
where the third equality is obtained by adding zero and the last equality holds because the Euler-characteristic of $\Gamma'$ yields $b^1(\Gamma')= E-k+1$.
\end{proof}

\begin{definition}
 Let $h:\Gamma \rightarrow \CL$ be a cover of degree $d$. The weights of the ends mapping to the three rays $u$, $v$ and $w$ of $\CL$ yield three partitions $\Delta_u$, $\Delta_v$ and $\Delta_w$ of $d$. The triple $\Delta=(\Delta_u,\Delta_v,\Delta_w)$ is called the \emph{ramification profile} of $h$.
\end{definition}

For example, the cover considered in \ref{exmp:cover} has ramification profile
$(\hspace{-0.0pt}(\hspace{-0.0pt}1,3\hspace{-0.0pt})\hspace{-0.0pt},\hspace{
-0.0pt} (\hspace{-0.0pt}4\hspace{-0.0pt})\hspace{-0.0pt},
\hspace{-0.0pt}(\hspace{-0.0pt}4\hspace{-0.0pt})\hspace{-0.0pt})$.

\begin{definition}
 We say that a combinatorial type $\alpha$ is \emph{trivalent}, if all vertices $V$ mapping to the center $c$ of $\CL$ have RH-number $r_V=0$ and all other vertices are trivalent and of genus zero (and thus have RH-number $r_V=1$).
\end{definition}

\begin{remark}\label{rem-wieners}
Automorphisms of a trivalent type arise only due to wieners (Figure \ref{fig-wieners}, see also  \cite{CJM10}).

% \begin{enumerate}
%   \item[(i)] Two inner edges of same weight, connecting the same two inner
% vertices
%  (we will call these objects \emph{wiener}, see also \cite{CJM10}):
\begin{figure}
 \begin{center}
 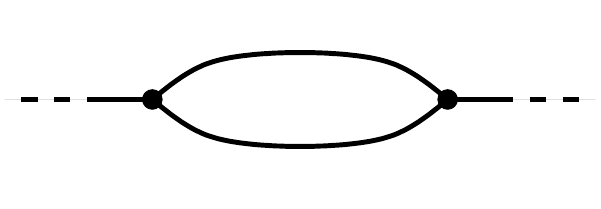
\end{center}
\caption{Wieners.}\label{fig-wieners}
\end{figure}
% \item[(ii)] Two ends of same weight adjacent to the same inner vertix, that is
% not
% mapped to $c$ (called \emph{balanced fork}, see also \cite{CJM10}:
% \begin{center}
%  \input{exmp:balancedfork.pdf_t}
% \end{center}
% \item[(iii)] Ends of same weight adjacent to the same vertex over $c$, that are
% mapped
% to the same end of $\CL$ (called \emph{balanced brush}):
% \begin{center}
%  \input{exmp:balancedbrush.pdf_t}
% \end{center}
%  \end{enumerate}
\end{remark}

We now take the set of all cones $D_\alpha$ such that the combinatorial type $\alpha$ is trivalent or a contraction of a trivalent type.
We glue these cones by identifying points on the boundary of $\overline{D}_\alpha$ with the corresponding point in the cone of its contraction as in Definition \ref{def-contraction}.
With this identification, the set of cones becomes an abstract polyhedral complex in the sense of \cite[Definition 3.4]{KM06} that we call $\Mg(\CL,\Delta)$, the
\emph{moduli space of tropical covers of $\CL$ of genus $g$ and with
ramification profile $\Delta$}.

By definition, the cones corresponding to trivalent types are the maximal cones of $\Mg(\CL,\Delta)$. 
We compute the dimension of a maximal cone:

\begin{lemma}\label{lem-dim}
 Let $\alpha$ be a trivalent type of degree $d$ genus $g$ covers of $\CL$ with
profile $\Delta$. Then the dimension of the cone $D_\alpha$ of $\alpha$ is
$\dim(D_\alpha) =     \#\Delta+2g-2-d$.
\end{lemma}

\begin{proof}
 For a trivalent cover of type $\alpha$, the dimension $\dim(D_\alpha)$ clearly equals the number of vertices which are not mapped to $c$: we can vary the lengths of the edges, staying within $D_\alpha$, in such a way that the images of these vertices move on $\CL$ (see also remark \ref{rem-numeq}). Each such moving image yields one degree of freedom. 
It follows that for a trivalent graph, the dimension $\dim(D_\alpha)$ equals the total number of vertex labels. The star-shaped cover --- i.e. the cover with one interior vertex adjacent to ends of weights $\Delta$ which are mapped to $u$, $v$ and $w$ accordingly --- is a contraction of every trivalent type. Since contraction by definition preserves the number of vertex labels, we can compute the number of vertex labels of the star-shaped cover in order to obtain the number of vertex labels of any trivalent cover. By the RH-condition, the star-shaped cover has $\#\Delta + 2g-2-d$ labels. The claim follows.
\end{proof}
Below, we equip each maximal cone of $\Mg(\CL,\Delta)$ with a weight, so that we can conclude the following result about the structure of $\Mg(\CL,\Delta)$:

\begin{thm}
The moduli space $\Mg(\CL,\Delta)$ of tropical covers of $\CL$ of genus $g$ and with
ramification profile $\Delta$ is an abstract weighted polyhedral complex of pure dimension $ \#\Delta+2g-2-d$.
\end{thm}

To introduce the weights of maximal cones, we need the following preparations.

Let $f:\ZZ^n\rightarrow \ZZ^m$ be a linear map. We call the index of $f$, $I_f$, the index of the sublattice $f(\ZZ^n)$ inside $\ZZ^m$.
\begin{definition}\label{def-equations}
 Let $\alpha$ be a combinatorial type of cover. In the underlying graph $\Gamma$, identify all vertices mapping to $c$ to one vertex. We call the graph obtained in this way $\Gamma'$. Pick $b^1(\Gamma')$ independent cycles, i.e.\ generators of $H_1(\Gamma',\ZZ)$. Each such generator is given as a chain of directed edges around the loop. In order to obtain a cover of type $\alpha$, we can choose lengths for the bounded edges, but we cannot choose them independently. The condition can be rephrased by stating that the images of the loops of $\Gamma'$ have to close up. In this way, we obtain $b^1(\Gamma')$ independent integral linear equations that cut out $D_\alpha$ from $(\RR_{>0})^{B}$, where $B$ denotes the number of bounded edges. We use the integral equations as defined by the weights of the edges that appear and do not cancel common factors.
\end{definition}

\begin{remark}\label{rem-numeq}
 An Euler-characteristic computation for $\Gamma'$ minus its ends shows that $1-b^1(\Gamma')=1+\#\{V|h(V)\neq c\}-B$, i.e.\ the number of equations in Definition \ref{def-equations} equals $B-\#\{V|h(V)\neq c\}$ (here, $B$ denotes again the number of bounded edges). It follows that the dimension of $D_\alpha$ equals $B-(B-\#\{V|h(V)\neq c\})=\#\{V|h(V)\neq c\}$. Indeed, in $D_\alpha$ we can vary the images of the vertices not mapped to $c$, and we have used this fact in the proof of lemma \ref{lem-dim}.
\end{remark}

\begin{example}\label{exmp:equations}
 Consider the cover of example \ref{exmp:cover}. We enumerate the inner edges
as indicated in Figure \ref{fig-equations} by the black numbers and red numbers. 
The black numbers represent a set of edges forming a spanning tree of $\Gamma'$, i.e.\ a set of edges whose lengths we can vary independently. The red edges each close a loop in $\Gamma'$, i.e.\ they depend on the lenths of the black edges.
\begin{figure}
\begin{center}
 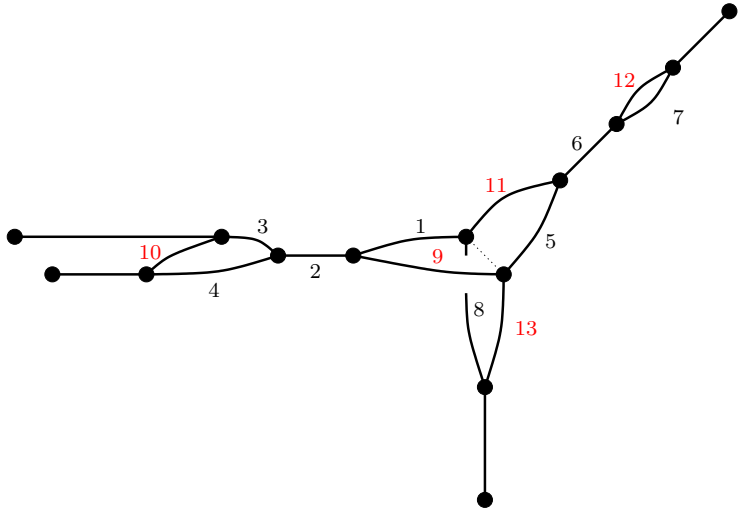
\end{center}\caption{Equations for cutting out $D_\alpha$.}\label{fig-equations}
\end{figure}
Denoting $x_{i}=l(e_i)$ we get the following six linear equations that
cut out $D_\alpha$ from $\RR_{>0}^{13}$:
\begin{eqnarray*}
 x_1-3\textcolor{red}{x_9}&=&0\\
 2x_3-2x_4+\textcolor{red}{x_{10}}&=&0\\
 3x_5-\textcolor{red}{x_{11}}&=&0\\
3x_7-\textcolor{red}{x_{12}}&=&0\\
x_8-3\textcolor{red}{x_{13}}&=&0.
\end{eqnarray*}
\end{example}

\begin{definition}\label{def-index}
 For a combinatorial type $\alpha$, we define $I_\alpha$ to be the index of the linear map $A_\alpha$ defined by the equations from definition \ref{def-equations}. 
\end{definition}

Note that while the matrix $A_\alpha$ depends on the choice of
generators of $H_1(\Gamma',\ZZ)$, its minors and therefore the index $I_\alpha$
do not (see \cite[chapter 5]{CJM10}).

% \begin{definition}\label{def-localautom}
% Let $h:\Gamma\rightarrow \CL$ be a cover and $V$ a vertex mapping to $c$ of local degree $d'$. Let $\Delta'_u$, $\Delta'_v$ and $\Delta'_w$ be the partitions of $d'$ corresponding to the weights of edges adjacent to $V$ mapping to $u$, $v$ and $w$ respectively. We define $\Aut_V(h)$ to be the \emph{local automorphisms} of $h$ at $V$, where $|\Aut_V(h)|=|\Aut(\Delta'_u)|\cdot|\Aut(\Delta'_v)|\cdot|\Aut(\Delta'_w)|$.
% 
% \end{definition}
% \begin{example}
% In example \ref{exmp:localCover} we have seen three local parts of tropical
% covers. Obviously the first two do not have local automorphisms. The last on the
% other hand yield the partition$(1,1)$ over the $u$-end and therefore
% $\Aut_V(h)=2$.
% \end{example}

Interestingly, we use (algebraic) Hurwitz numbers to define the weights of the maximal cones of $\Mg(\CL,\Delta)$. We first recall the definition of Hurwitz numbers:

\begin{definition}\label{def-algnumb}
 Fix $r$ points  $p_1, \ldots, p_r$ in $\mathbb{P}^1\setminus\{0,1,\infty\}$, and  $\Delta=(\Delta_u,\Delta_v,\Delta_w)$, where $\Delta_u$, $\Delta_v$ and $\Delta_w$ are partitions of the same integer $d$. Then the \textit{(triple) Hurwitz number} $H^{g}_d(\Delta)$ is defined as the weighted number of degree $d$ covers $f: C\rightarrow \mathbb{P}^1$ satisfying:
\begin{itemize}
 \item C is a marked smooth connected curve of genus $g$ (all preimages of $0$, $1$ and $\infty$ are marked);
\item $f$ ramifies with profile $\Delta_u$ over $0$;
\item $f$ ramifies with profile $\Delta_v$ over $1$;
 \item $f$ ramifies with profile $\Delta_w$ over $\infty$;
\item $f$ has simple ramification over $p_i$ for all $i=1,\ldots,r$.
\item $f$ is unramified over $\mathbb{P}^1\setminus\{p_1,\ldots, p_r,0,1,\infty\}$;
\end{itemize}
Each cover is weighted by $\frac{1}{\Aut(f)}$.
\end{definition}
Note that we mark the preimages of the three special ramification points. In the
literature, the analogous definition where the preimages are not marked is also
common. It differs from our definition just by a factor of
$\frac{1}{\Aut(\Delta)}$.

The (algebraic) Riemann-Hurwitz formula (see e.g.\ \cite{Har77}, Corollary IV.2.4) states that $2g-2-d+ \#\Delta=r$.

\begin{remark}\label{rem-symmetricgroup}
 By matching a cover with a monodromy representation, we can count Hurwitz numbers by counting tuples of elements of the symmetric group $\mathbb{S}_d$ that multiply to the identity and act transitively on $\mathbb{S}_d$, see e.g.\ \cite{v:tmsocagwt}, Proposition 3.17, or \cite{Joh12}. The transitivity condition is
equivalent to connectedness of the source curve.
\end{remark}

\begin{definition}
 Let $h:\Gamma\rightarrow \CL$ be a trivalent cover and $V$ a vertex mapping to $c$ of local degree $d'$. As before, let $\Delta'_u$, $\Delta'_v$ and $\Delta'_w$ be the partitions of $d'$ corresponding to the weights of edges adjacent to $V$ mapping to $u$, $v$ and $w$ respectively. Since the RH-number of $V$ is zero, we have 
$\#\Delta'_u+\#\Delta'_v+\#\Delta'_w+2g_V-2-d'=0$.

We define $H_V:=H^{g_V}_{d'}(\Delta'_u,\Delta'_v,\Delta'_w)$ to be the (algebraic) triple Hurwitz number of genus $g_V$ degree $d'$ covers of $\PP^1$ with ramification profile $\Delta'_u$ over $0$, $\Delta'_v$ over $1$ and $\Delta'_w$ over $\infty$ (see Definition \ref{def-algnumb}). The (algebraic) Riemann-Hurwitz formula implies $2g_V-2-d'+ \#\Delta'_u+\#\Delta'_v+\#\Delta'_w-s=0$, where $s$ denotes the number of simple ramifications of a cover of degree $d'$, genus $g_V$ and with three special ramification points of profiles $\Delta'_u$, $\Delta'_v$ and $\Delta'_w$. Hence the above implies $s=0$, i.e.\ a cover with such three special ramification profiles has no other ramification.
\end{definition}

\begin{definition}\label{def-weights}
For a maximal cone of $\Mg(\CL,\Delta)$, resp.\ for a trivalent type $\alpha$, we define its weight $\omega(\alpha)$ to be 
$$\omega(\alpha):= \frac{1}{2^k}\cdot I_\alpha \cdot \prod_V  H_V,$$
where $k$ denotes the number of wieners (see remark \ref{rem-wieners}) and the product runs over all vertices $V$ mapping to $c$. 
\end{definition}

Note that this definition is natural when compared to other definitions of weights in tropical moduli spaces, see e.g.\ \cite[definition 3.5]{KM06} or \cite[definiton 5.10]{CJM10}. Also, it is natural from the point of view of tropical intersection theory, since the cones $D_\alpha$ are cut out by the equations of index $I_\alpha$.

\subsection{Comparing the weights}
In \cite{GMO}, moduli spaces of rational tropical covers of a line are constructed as balanced polyhedral complexes in a surrounding vector space. The approach taken there is different from ours, although also motivated by tropical intersection theory. The moduli spaces are constructed by gluing pieces using tropical intersection theory, respectively by interpreting them as a marked polyhedral subcomplex of the moduli space of rational curves. The methods used in \cite{GMO} highly rely on the genus to be zero and cannot be easily generalized to arbitrary genus. 
For the case of genus zero, we now show that the weight of a cone $D_\alpha$ that we define in definition \ref{def-weights} coincides with the weight in \cite{GMO} (which is shown there to satisfy the balancing condition). This may serve as an additional justification for our choice of weights.
In \cite{GMO}, a cone $D_\alpha$ of a combinatorial type $\alpha$ obtains the weight
\begin{equation}\prod_V H_V\cdot \gcd_{T}(\prod_{e\in T^c}w(e)),\label{eq-weight}
\end{equation}
where the $\gcd$ is taken over all spanning trees $T$ of the graph $\Gamma'$ obtained as in definition  \ref{def-equations} from the underlying graph $\Gamma$ and the product goes over all vertices $V$ which are mapped to $c$.

\begin{prop}
For a combinatorial type $\alpha$ of a rational cover of $\CL$, the weight we define for $D_\alpha$ in definition \ref{def-weights} coincides with the weight of equation \ref{eq-weight} which is shown in \cite{GMO} to satisfy the balancing condition.
\end{prop}
\begin{proof}
Note that since $\Gamma$ is rational there are no wieners, so we only need to show that $I_\alpha$ equals 
$\gcd_{T}(\prod_{e\in T^c}w(e))$.
Let $\Gamma'$ and $b$ be as in definition \ref{def-equations}, and
$g'=\#\{b^1(\Gamma')\}$ the number of equations. Moreover we label the columns
of $A_\alpha$ (see definition \ref{def-index}) by the corresponding edge in $\Gamma'$. As in
\cite[Lemma 3.20]{GS11} $I_\alpha$ equals the greatest common divisor of
its $g'\times g'$-minors. Thus we have to show $$\gcd_{1\leq
i_1<\ldots<i_{g'}\leq b}M_i=\gcd_{T}(\prod_{e\in T^c}w(e)),$$ where
$M_i=M_{(i_1,\ldots,i_{g'})}$ is the minor of $A_\alpha$ containing the columns
$e_{i_1},\ldots,e_{i_{g'}}$ and $T$ goes over all spanning trees of
$\Gamma'$. 

We show that the set of nonzero minors $M_{(i_1,\ldots,i_{g'})}$ coincides with the set of all $ \prod_{e\in T^c}w(e)$ for spanning trees $T$ of $\Gamma'$.

Let $M_{(i_1,\ldots,i_{g'})}$ be any minor and denote
$S=\{e_{i_1},\ldots,e_{i_{g'}}\}$.

Assume first $S^c$ is not a tree in
$\Gamma'$. Since $\Gamma'$ is a connected graph of genus $g'$ and we remove $g'$ edges to obtain $S^c$, the assumption implies that $S^c$ contains a cycle. Since the minors do
not depend on the choice of equations, we can choose the first row to correspond to the cycle in $S^c$.
In this row, we then have zeros in the columns corresponding to
$S$, so $M_i$ is zero.

Now let $T:=S^c$ be a tree. Adding an edge $e$ of $T^c$ to $T$ produces
a unique cycle containing $e$ and no other edge of $T^c$. 
We use these $g'$ fundamental cycles of $T$ to write down the matrix $A_\alpha$.
We obtain
$$A_\alpha=\left(\begin{array}{ccc|ccc}
             w(e_{i_1})&&0\\
	    &\ddots&&\hspace{10pt}&\ast&\hspace{10pt}\\
	    0&&w(e_{i_{g'}})
            \end{array}\right),$$
where the first $g'$ columns correspond to the edges $e_{i_1},\ldots,e_{i_{g'}}$ in $S=T^c$
and the remaining columns to edges in $T$. Then clearly
$M_i=\prod_{k=1}^{g'}w(e_{i_k})=\prod_{e\in T^c}w(e)$.
Thus every nonzero minor equals $\prod_{e\in T^c}w(e)$ for some spanning tree $T$. Vice versa, the above construction also shows that for a given tree $T$ we obtain a minor which equals $\prod_{e\in T^c}w(e)$.
\end{proof}

\section{The branch map}

\begin{definition}
 The \emph{(tropical) branch map} on the moduli space $\Mg(\CL,\Delta)$ is
defined as
\begin{eqnarray*}
 \br^{\trop}:\Mg(\CL,\Delta) & \rightarrow & \CL^r\\
(h:\Gamma\rightarrow \CL) & \mapsto & (h(V_1),h(V_2),\ldots,h(V_r)),
\end{eqnarray*}
where $r=\#\Delta+2g-2-d$ is the total number of labels (and the
dimension of $\Mg(\CL,\Delta)$).
\end{definition}
It follows easily that $\br^{\trop}$ is a morphism of weighted polyhedral complexes of the same dimension in the sense of \cite[Definition 4.1]{KM06}.

Remember that the degree of a morphism $f$ of weighted polyhedral complexes of the same dimension is defined to be the sum of the weights of cones times the local multiplicities of cones (we denote the latter by $\mult_D f$ for a cone $D$), where the sum goes over all inverse images of a point in general position \cite[Definition 4.1]{KM06}, i.e.\
$$ \deg(f)= \sum_{Q\;|\; f(Q)=P} \omega(D(Q)) \mult_{D(Q)} f,$$ where $D(Q)$ denotes the maximal cone that contains $Q$ in its interior.

\begin{lemma}\label{lem-productedges}
 Let $\alpha$ be the type of a trivalent cover. We have
$$I_\alpha\cdot \mult_{D_\alpha}\br^{\trop}=\prod_{e}w_e,$$
where the product goes over all bounded edges $e$ of the underlying graph
$\Gamma$, $w_e$ denotes their weights, $I_\alpha$ the lattice index defined in
\ref{def-index} and $\mult_{D_\alpha}\br^{\trop} $ the local multiplicity of the branch
map just as above.
\end{lemma}
\begin{proof}
 This is a straight-forward generalization of Remark 5.19 and Lemma 5.26
of \cite{CJM10}. 
\end{proof}

 We now state our main result.
\begin{thm}\label{thm}
The degree of $\br^{\trop}$ is constant, i.e.\ it does not depend on the choice of the point in general position that we pull back.
\end{thm}

As a consequence of Theorem \ref{thm}, we can define:

\begin{definition}
 For $g$ and $\Delta$, we define the tropical Hurwitz number $H^{g,\trop}_d(\Delta)$ (where $d$ is the sum of the parts of $\Delta_u$, $\Delta_v$ and $\Delta_w$, resp.) to be the degree of the branch map $\br^{\trop}:\Mg(\CL,\Delta)  \rightarrow  \CL^r$.
\end{definition}

\begin{lemma}\label{lem-BBM}
 Our definition of tropical Hurwitz number agrees with the definition of
\cite{BBM10}, up to a factor of $|\Aut(\Delta)|$ that arises
because we mark the ends.
\end{lemma}

\begin{proof}
The definition of tropical Hurwitz number in \cite{BBM10} counts covers where all the ramification data is imposed at the ends, i.e.\ simple ramification in the interior appearing as trivalent vertices is not considered. To interpret our covers in this context, we need to add an extra end to
$\CL$ at the image of every trivalent vertex $V$ not mapping $c$, and
analogously add one (unmarked) end of weight $2$ and $\deg_hV-2$  (unmarked) ends of weight
$1$ to  $V$ as well as $\deg_hW$  (unmarked) ends of weight one to every $W\in
h^{-1}(h(V))\setminus\{V\}$. We call the new tropical curves obtained in this way $\Gamma'$ and $\CL'$
respectively. We extend $h$ to a cover $h':\Gamma'\rightarrow\CL'$ such that the
new ends of $\Gamma'$ are mapped to the new ends of $\CL'$ in the obvious manner.

By Lemma \ref{lem-productedges} we can write the contribution of each combinatorial type
of cover to our count as
\begin{equation}\omega(\alpha)\cdot \mult_{D_\alpha}\br^{\trop}= \frac{1}{2^k}\cdot \prod_e
w_e \cdot \prod_V  H_V,\label{eq1}\end{equation}
where $k$ denotes the number of wieners. Note that by remark \ref{rem-wieners} $|\Aut(h)|=2^k$.
On the other hand, $h':\Gamma'\rightarrow\CL'$ is counted in \cite{BBM10} with multiplicity
\begin{equation}\frac{1}{|\Aut(h')|}\cdot \prod_{e'}
w_{e'} \cdot \prod_{V'} H_{V'}.\label{eq2}\end{equation}
(Note that in \cite{BBM10} the authors work with a definition of (algebraic) Hurwitz numbers where we do not mark the preimages of the three special ramification points, consequently they have to multiply their Hurwitz number with a factor reflecting the local automorphisms, i.e.\ the automorphisms of the three local partitions.)
The automorphisms $\Aut(h')$ here consist of automorphisms of the unmarked ends, and the wieners as before.
We now analyze the difference between the two expressions.

\begin{itemize}
 \item Assume $V$ is a trivalent vertex not mapping to $c$, then $V$ does not
contribute any Hurwitz number to (\ref{eq1}). In (\ref{eq2}) the corresponding vertex $V'$ provides a factor $  H_{V'}= (\deg_h(V)-1)!$ (this number reflects the number of ways to mark the preimages of the simple branch point). This factor is annihilated by the corresponding global automorphisms in the whole product. 
%The latter
%is $2$ if and only if $V$ is adjacent to two edges with identical weight
%(locally a balanced fork). But then $H_{V'}=\frac{1}{2}$ which cancels the contribution from the automorphisms. If $h$ has no local automorphisms at $V$ both contributions are $1$.
\item Let $W\in h^{-1}(h(V))\setminus\{V\}$ for a vertex $V$ as above. Similar
to the former case $W$ yields no contribution to (\ref{eq1}) and in (\ref{eq2}) we get $H_{W'}=\frac{1}{\deg_hW}\cdot (\deg_hW-1)!$. By adding the extra ends at $W$ we subdivide an edge $e$ of $\Gamma$ into two edges providing an additional factor of $w_e= \deg_hW$ which together with the new global automorphisms cancels the contribution of $H_{W'}$.
% Similar
% to the former case $W$ yields no contribution to (\ref{eq1}) and in (\ref{eq2}) the local automorphisms of $h'$ at the corresponding vertex $W'$ annihilate the contribution to the global automorphisms coming from this vertex. The number $H_{W'}$ equals $\frac{1}{\deg_hW}$. By adding the extra ends at $W$ we subdivide an edge $e$ of $\Gamma$ into two edges providing an additional factor of $w_e= \deg_hW$ which cancels the contribution of $H_{W'}$.
\end{itemize}

Furthermore, vertices mapping to $c$ yield the same contributions to both counts (\ref{eq1}) and (\ref{eq2}).
We count covers with marked ends, so for each cover $h'$ we have to multiply by a factor taking into account the possibilities to mark the ends. This factor times the contribution to $|\Aut(h')|$ arising from these ends equals $|\Aut(\Delta)|$. The contribution to $|\Aut(h')|$ of the newly attached ends all cancel as discussed above. There remain only contributions from wieners, which we also have in (\ref{eq1}).
 It follows that the two expressions agree up to a factor of $|\Aut(\Delta)|$.
\end{proof}

As a consequence of Lemma \ref{lem-BBM}, we can conlude:
\begin{thm}
The tropical Hurwitz numbers $H^{g,\trop}_d(\Delta)$ that we define using tropical intersection theory on an appropriate moduli space equal their algebraic counterparts $H^{g}_d(\Delta)$ (see Definition \ref{def-algnumb}).
\end{thm}
\begin{proof}
 This follows from the Correspondence Theorem 2.11 in \cite{BBM10}.
\end{proof}

The main ingredient for the proof of Theorem \ref{thm} is a duality between tropical resolutions of a codimension-one-case and boundary points of a one-dimensional algebraic moduli space. 
We first explain the one-dimensional case in detail before deducing the consequences for the general situation.
\subsection{The one-dimensional case}\label{one-dimCase}
Throughout this subsection, fix a
ramification profile $\Delta$, a degree $d$ and a genus $g$ such that
$\#\Delta+2g-d-2=1$, i.e. the covers in $\Mg(\CL,\Delta)$ have exactly one
label. Then there is exactly one (combinatorial type of) cover that is not
trivalent, namely the star-shaped cover with a vertex of genus $g$ and $1$
label over $c$. Obviously, $\Mg(\CL,\Delta)$ as abstract polyhedral complex is just a star itself: a collection of one-dimensional rays adjacent to the star-shaped curve. Each ray corresponds to a possible resolution of the star-shaped curve, i.e.\ to a cover of $\CL$ with one trivalent vertex mapping to one of the rays of $\CL$.
Topologically, there are three different types for such resolutions:
we can either join two edges (as e.g.\ in the top row
on the left of Figure \ref{fig:tropResolutions}), split
an edge while extracting genus from the vertex over $c$ (as e.g.\ in the bottom
row on the left of Figure \ref{fig:tropResolutions}) or split an edge
and the interior vertex (as e.g.\ in the top row in the middle of Figure
\ref{fig:tropResolutions}).

\begin{example}\label{example:tropicalResolutions}
Consider the space $\Mone(\CL,((3,1,1),(5),(3,2)))$. The star-shaped
combinatorial type in this space has an interior vertex of genus one. Its
resolutions, i.e.\ the
trivalent combinatorial types in this space --- ordered by the position of the
image $p$ of their labelled point on the different ends of $\CL$ --- are
depicted in Figure \ref{fig:tropResolutions}. (As before blue numbers denote
edge weights and red numbers are
the genus on vertices.)
\begin{figure}
\begin{center}
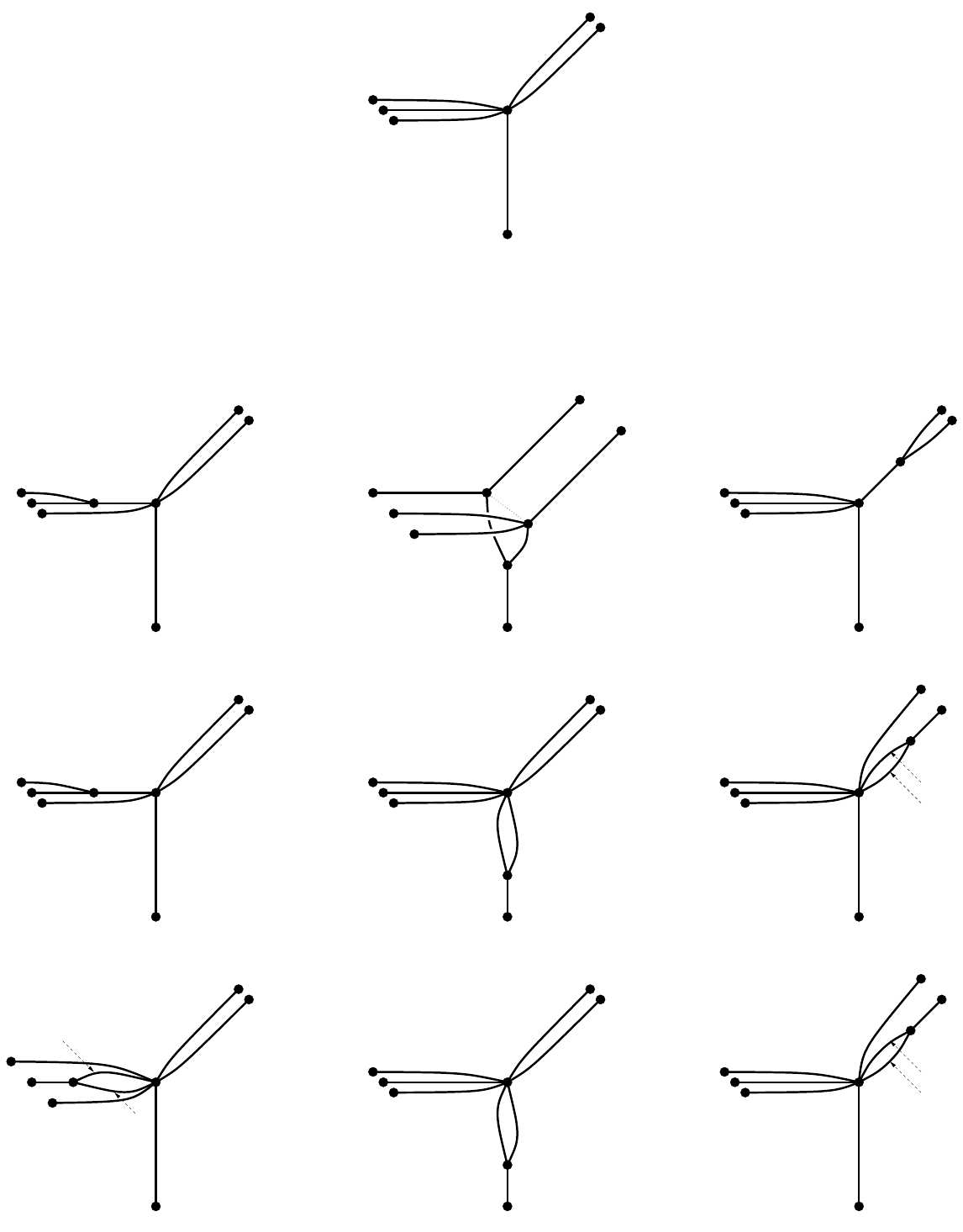
\end{center}
\caption{Resolutions of the star-shaped cover in $\Mone(\CL,((3,1,1),(5),(3,2)))$.}
\label{fig:tropResolutions}
\end{figure}
In the picture, we neglect the marking of the ends as usual. This implies that
e.g. the picture in the top row on the left actually combines two marked
pictures, for the two possibilities to mark the two ends of weight one.
\end{example}

\begin{definition}
Fix the three points $0$, $1$ and $\infty$ in $\PP^1$. We consider relative stable maps to $\PP^1$, relative to these three points, with profiles $\Delta_u$, $\Delta_v$ and $\Delta_w$ respectively. We denote the space of such relative stable maps by $\overline{M}_{g,\#\Delta, \Delta}(\PP^1,d)$.
\end{definition}
\begin{remark}\label{rem-algspace}
The space $\overline{M}_{g,\#\Delta, \Delta}(\PP^1,d)$ is a one-dimensional moduli stack (\cite{Li01}, \cite{Li02}, for a nice introduction to relative stable maps to $\PP^1$ see also \cite{v:tmsocagwt}).
Points in $\overline{M}_{g,\#\Delta, \Delta}(\PP^1,d)$ roughly correspond to maps of a source $C$ to a chain of $\PP^1$s such that the kissing condition is satisfied above each node, the three ramification profiles are satisfied about two points in the first copy of $\PP^1$ of the chain and one point in the last copy and all preimages of the branch points are marked. The stability condition implies that there is at most one node in the target, i.e.\ at most two copies of $\PP^1$. 
Points in the interior are ramified covers of $\PP^1$ with the three special ramification profiles as above and one further simple ramification at a point $t\neq \{0,1,\infty\}$. 
At the boundary, i.e.\ when $t$ moves to one of the three points, the covers degenerate to covers of two copies of $\PP^1$ as follows. Consider the situation where $t$ moves to $0$. 
Then we have covers of a chain of two $\PP^1$s that satisfy the kissing condition above the node, say the ramification profile above the node is $\tilde \Delta$. On one copy of $\PP^1$, we then have three ramification points with profiles $\Delta_u$, a simple ramification and $\tilde \Delta$. On the other, we have $\tilde \Delta$, $\Delta _v$ and $\Delta_w$.
The possibilities for $\tilde \Delta$ are restricted by the cut-and-join
relations: to obtain $\tilde\Delta$, we can either divide one entry of
$\Delta_u$ into two parts or sum two parts of $\Delta_u$. This follows from
remark \ref{rem-symmetricgroup}: by matching a cover with a tuple of elements in
the symmetric group, the simple ramification corresponds to a transposition
$\tau$ while $\Delta_u$ and $\tilde\Delta$ correspond to permutations $\sigma_u$
and $\tilde\sigma$ of appropriate cycle type satisfying
$\sigma_u\circ\tau=\tilde\sigma$. A transposition can either cut a cycle or join
two cycles in a permutation.
\end{remark}
The duality between boundary points of $\overline{M}_{g,\#\Delta, \Delta}(\PP^1,d)$ and rays of $\Mg(\CL,\Delta)$ goes by the dual graph construction:

\begin{definition} For an element of $\overline{M}_{g,\#\Delta, \Delta}(\PP^1,d)$, we construct its dual graph as follows:
\begin{itemize}
 \item For every component $C_i$ of the source curve $C$, we draw a vertex with genus $g(C_i)$;
\item for every node of component $C_i$ and $C_j$ we draw an edge between the vertices $i$ and $j$, the weight of the edge equals the intersection multiplicity of the two components at the node;
\item for every marked point on $C_i$ we draw a marked end adjacent to the vertex $i$, the weight of the end equals the ramification index at the marked point.
\end{itemize}
We straighten two-valent vertices. We interpret the outcome as a combinatorial type of tropical covers of $\CL$ by mapping the parts that go to $0$, $1$ and $\infty$ to $u$, $v$ and $w$ respectively.
\end{definition}
Obviously, the dual graph of a cover in the interior of $\overline{M}_{g,\#\Delta, \Delta}(\PP^1,d)$ is just the star-shaped cover in $\Mg(\CL,\Delta)$.

\begin{example}
Figure \ref{fig:classicalCover} sketches a cover corresponding to an interior point of 
$\overline{M}_{1,\#\Delta, \Delta}(\PP^1,5)$ where
$\Delta=((3,1,1),(5),(3,2))$. A cover corresponding to a boundary point in the
moduli space is sketched in figure \ref{fig:classicalResolution}. We neglect markings as usually.
\begin{figure}
 \begin{center}
 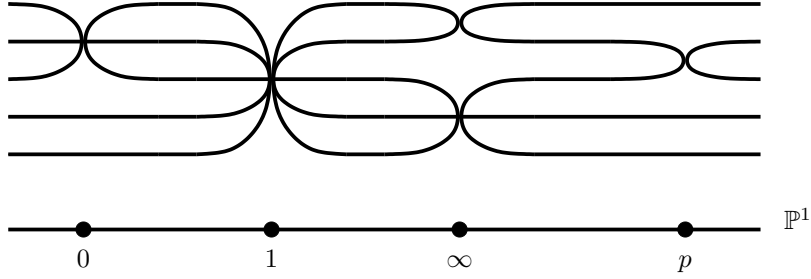
 \caption{An algebraic cover dual to the star-shaped cover in
figure \ref{fig:tropResolutions}.}\label{fig:classicalCover}
\end{center}
\end{figure}
\begin{figure}
 \begin{center}
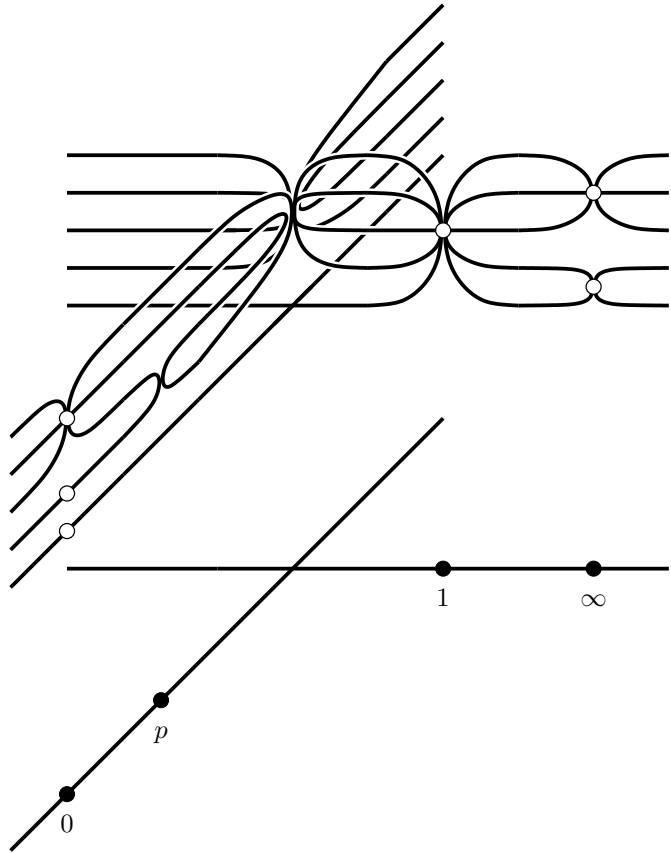
 \caption{The boundary point of $\overline{M}_{1,\#\Delta, \Delta}(\PP^1,5)$ dual to the tropical cover on the left of the first
row in figure \ref{fig:tropResolutions}.}\label{fig:classicalResolution}
\end{center}
\end{figure}

Figure \ref{fig:classicalResolutions} very roughly sketches all covers
corresponding to boundary points of $\overline{M}_{1,\#\Delta, \Delta}(\PP^1,5)$
dual to the tropical covers in figure \ref{fig:tropResolutions}. The order is
the same in both pictures. Also here, we neglect the markings of the preimages
of the three special branch points, the picture on the top left actually
combines two marked pictures. The top left picture represents the same cover as
figure
\ref{fig:classicalResolution}, the kissing condition is indicated by the broken
line.
 
\begin{figure} 
  \begin{center}   
   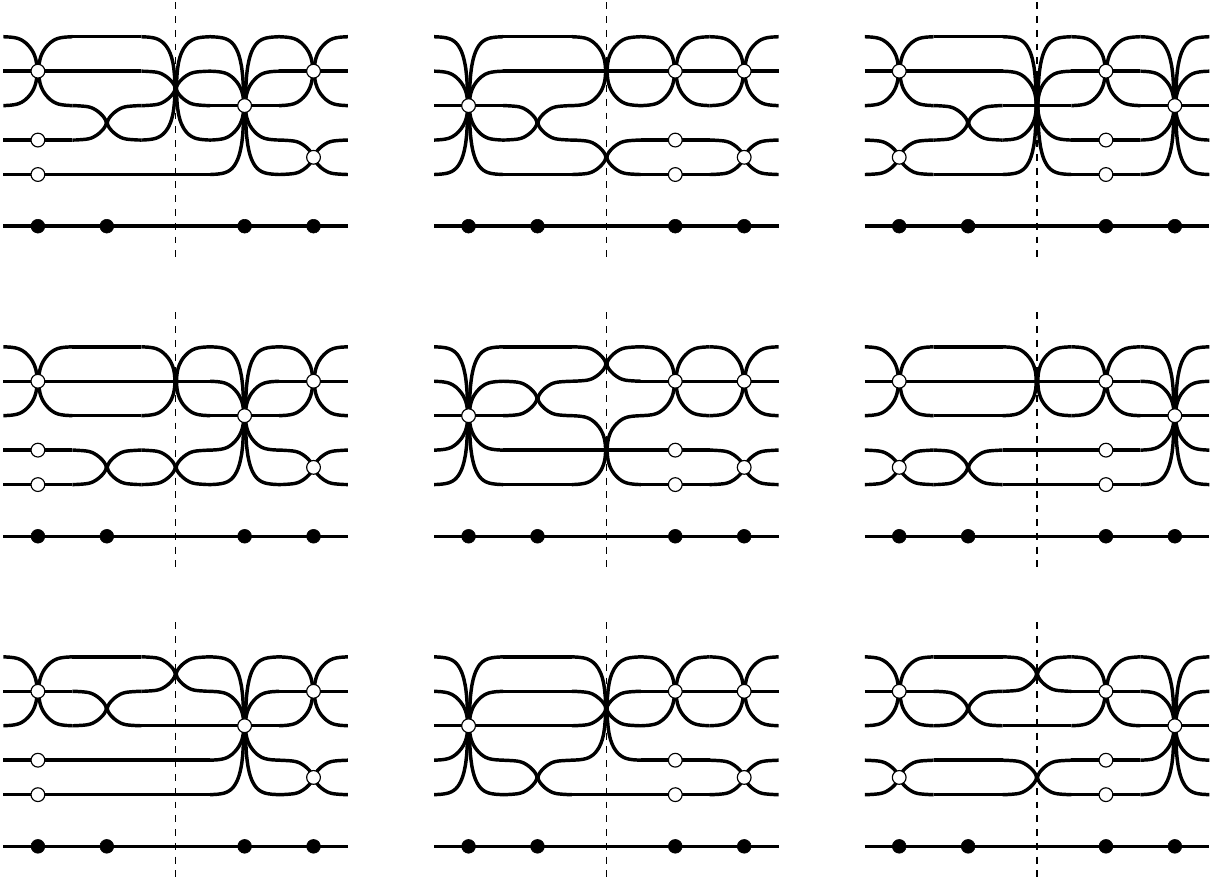
  \end{center}
\caption{Boundary points of $\overline{M}_{1,\#\Delta, \Delta}(\PP^1,5)$ dual to the tropical
resolutions in figure \ref{fig:tropResolutions}.}\label{fig:classicalResolutions}
 \end{figure}
\end{example}

\begin{prop}\label{prop:1-1-correspondence}
 The boundary points of $\overline{M}_{g,\#\Delta, \Delta}(\PP^1,d)$ are in
$1:1$-corres-pondence with rays of $\Mg(\CL,\Delta)$ via the dual graph construction.
More precisely, boundary points where $t$ goes to $0$ correspond to tropical
covers with a trivalent vertex above $u$, where $t$ goes to $1$ to covers with a
trivalent vertex above $v$ and where $t$ goes to $\infty$ to covers with a
trivalent vertex above $w$.
\end{prop}

\begin{proof}
 Take a point in the boundary of $\overline{M}_{g,\#\Delta, \Delta}(\PP^1,d)$, say where $t$ moved to $0$.
We claim that the dual graph $\Gamma$ is a possible resolution of the star-shaped cover in $\Mg(\CL,\Delta)$ with a trivalent vertex above $u$. As described in remark \ref{rem-algspace}, such a boundary point is a cover of two copies of $\PP^1$, one copy with ramification profiles $\Delta_u$, simple and $\tilde \Delta$, the other with $\tilde \Delta$, $\Delta_v$ and $\Delta_w$. The possibilities for $\tilde \Delta$ are given by the cut-and-join relations (see remark \ref{rem-algspace}) . A cover with profiles $\Delta_u$, simple and $\tilde \Delta$ contains one rational component $C_1$ with the simple ramification and two more ramification profiles, one totally ramified and the other in two parts. The dual vertex is a trivalent vertex of genus zero which is mapped to $u$. The remaining components are mapped trivially (and thus also rational), thus their dual vertex is two-valent, with one adjacent marked end and one bounded edge of the same weight connecting it to a vertex corresponding to a component 
covering the other copy of $\PP^1$. 
We have the following possibilities: 
\begin{enumerate}
 \item There is exactly one component covering the other copy of $\PP^1$, and it meets $C_1$ in two nodes. The dual graph then is as e.g.\ in the bottom
row on the left of Figure \ref{fig:tropResolutions}. 
\item There is exactly one component covering the other copy of $\PP^1$, and it meets $C_1$ in one node. Then $\tilde \Delta$ is obtained from $\Delta_u$ by summing two parts, and consequently we have two marked points in $C_1$. The dual graph is as e.g.\ in the top row
on the left of Figure \ref{fig:tropResolutions}.
\item There are two components covering the other copy of $\PP^1$, each meeting $C_1$ in one node. The dual graph is as e.g.\ in the top row in the middle of Figure
\ref{fig:tropResolutions}.
\end{enumerate}
Vice versa, we can obviously construct for each combinatorial type of tropical cover corresponding to a ray of $\Mg(\CL,\Delta)$ a boundary point in $\overline{M}_{g,\#\Delta, \Delta}(\PP^1,d)$ whose dual graph equals the combinatorial type.
\end{proof}

\begin{prop}\label{prop:brIsCover}
 The branch map $\br: \overline{M}_{g,\#\Delta, \Delta}(\PP^1,d) \rightarrow \PP^1$ taking a cover to the image of its simple branch point is itself a cover of $\PP^1$ of degree $H^g_d(\Delta)$, branched above $0$, $1$ and $\infty$. 
In particular, $\br^\ast(0)=\br^\ast(1)=\br^\ast(\infty)$, and each consists of the boundary points desribed above.
\end{prop}
\begin{proof}
 By \cite{fp:smabd}, $\br$ is a natural map of stacks (see also \cite{v:tmsocagwt}, Section 6.2). The statement about the degree and the branching is obvious.
\end{proof}

\begin{lemma}\label{lem:1:1corrWithMult}
 The multiplicity of a boundary point of $\overline{M}_{g,\#\Delta, \Delta}(\PP^1,d)$ in $\br^\ast(p)$ for $p=0$, $1$ or $\infty$ equals the tropical multiplicity of the combinatorial type of cover given by the dual graph.
\end{lemma}
\begin{proof}
We formulate the argument for $p=0$ to keep notation simple. Since $\br$ is a
branched cover, we can determine the multiplicity of a boundary point in
$\br^\ast(0)$ by counting the number of covers in $\overline{M}_{g,\#\Delta,
\Delta}(\PP^1,d)$ with the simple ramification at $t$ close to $0$ that
degenerate to the given boundary point. 

We count these covers in terms of monodromy representations as in remark \ref{rem-symmetricgroup}.  As every ramification point over $0$, $1$ and
$\infty$ is marked, we can think of $\Delta_u$, $\Delta_v$ and
$\Delta_w$ as marked partitions, where the marking is induced by the marks of the preimages of $0$, $1$ and $\infty$, respectively.
We also consider  permutations $\sigma\in\Perm_d$ together with a marking of their cycles and call this a \emph{marked permutation}. By abuse of notation, we still denote a marked permutation by $\sigma\in \Perm_d$.

We say that a marked permutation $\sigma$
is of \emph{marked cycle type} $\Delta_u$ and write
$\sigma\in\Perm_d^{(\Delta_u)}$ if the marked tuple of its cycle lengths agrees with
the marked partition $\Delta_u$.

Following remark \ref{rem-symmetricgroup}, the Hurwitz number $H^{g}_d(\Delta)$ equals

$$H^{g}_d(\Delta)=\frac{1}{d!}\cdot\#\left\{(\sigma_u,\sigma_v,\sigma_w,\tau)\right\},$$
where the tuples in the braces satisfy
\begin{itemize}
 \item $\sigma_u$, $\sigma_v$ and $\sigma_w$ are marked permutations satisfying $\sigma_u\in\Perm_d^{(\Delta_u)}$, $\sigma_v\in\Perm_d^{(\Delta_v)}$ and
$\sigma_w\in\Perm_d^{(\Delta_w)}$ respectively,
\item $\tau$ is an unmarked
transposition in $\Perm_d$,
\item $\tau\circ\sigma_u\circ\sigma_v\circ\sigma_w=\operatorname{id}
_{\Perm_d}$ and
\item $\langle\tau,\sigma_u,\sigma_v,\sigma_w\rangle_{\Perm_d}$ acts transitively on $\{1,\ldots,d\}$.
\end{itemize}

Now consider a possible kissing condition $\tilde\Delta$. As in remark \ref{rem-algspace}, it is obtained from $\Delta_u$ by either splitting one part into two or summing to parts to one. In the first case, we consider $\tilde\Delta$ as a partially marked partition (where the two new parts are not marked). Analogously, we also consider partially marked permutations and say they are of partially marked cycle type $\tilde\Delta$, if the partially marked partition of cycle lengths agrees with $\tilde\Delta$. By abuse of notation, we also write $\tilde\sigma\in\Perm_d^{(\tilde\Delta)}$ if $\tilde\sigma$ is of partially marked cycle type $\tilde\Delta$. In the following, it should always be clear from the context whether a permutation is marked, partially marked or unmarked.

Fix a boundary point in $\br^\ast(0)$ with kissing condition $\tilde\Delta$.
Remember from remark \ref{rem-algspace} that for a boundary point, the target consists of two copies of $\PP^1$ meeting in a node. One copy is covered with ramification profiles $\Delta_u$, simple and $\tilde \Delta$, the other by $\tilde \Delta$, $\Delta_v$ and $\Delta_w$. There is one component called $C_1$ above the first copy of $\PP^1$ which contains the simple ramification.

\hspace{0.2cm}

Assume first the dual graph of the boundary point is as in case (1) of the proof of proposition
\ref{prop:1-1-correspondence}, i.e.\ as e.g.\ in the bottom
row on the left of Figure \ref{fig:tropResolutions}. Then $\tilde\Delta$ is
obtained from $\Delta_u$ by splitting the part $m$ into positive intergers $m_1$
and $m_2$ with $m_1+m_2=m$.
  If we consider covers with simple ramification at $t$ close to $0$, we
can count the ones which degenerate to this boundary point as follows:

 \begin{equation}\frac{1}{d!}\cdot\#\left\{(\sigma_u,\sigma_v,\sigma_w,\tau)\left|
	    \begin{array}{l}
\bullet\ \sigma_u\in\Perm_d^{(\Delta_u)},\sigma_v\in\Perm_d^{(\Delta_v)},
\sigma_w\in\Perm_d^{(\Delta_w)}\\
\bullet\ \tau\textnormal{ an unmarked transposition in }\Perm_d\\
\bullet\
\tau\circ\sigma_u\circ\sigma_v\circ\sigma_w=\operatorname{id}_{\Perm_d}\\
\bullet\ \langle\tau,\sigma_u,\sigma_v,\sigma_w\rangle_{\Perm_d}\textnormal{
acts transitively on }\{1,\ldots,d\}\\
\bullet\ \langle\tau\circ\sigma_u,\sigma_v,\sigma_w\rangle_{\Perm_d}\textnormal{
acts transitively on }\{1,\ldots,d\}\\
\bullet\ \tau\circ\sigma_u\in\Perm_d^{(\tilde\Delta)}
             \end{array}
\right.\right\}.\label{eq-count}\end{equation}
The second transitivity condition reflects the fact there is only one component above the other copy of $\PP^1$ which meets $C_1$ in two nodes.
Obviously the first transitivity condition is obsolete. We can order the set
of tuples by the result of $\sigma_u\circ\tau$ and accordingly write the number as $\frac{1}{d!}$ times 
the sum over all $\tilde\sigma\in\Perm_d^{(\tilde\Delta)}$ of products of two factors:

$$\#\left\{
(\sigma_v,\sigma_w)\left|
	    \begin{array}{l}
\bullet\ \sigma_v\in\Perm_d^{(\Delta_v)},\sigma_w\in\Perm_d^{(\Delta_w)}\\
\bullet\ \sigma_w\circ\sigma_v\circ\tilde\sigma=\operatorname{id}_{\Perm_d}\\
\bullet\ \langle\tilde\sigma,\sigma_v,\sigma_w\rangle_{\Perm_d}\textnormal{
acts transitively on }\{1,\ldots,d\}
             \end{array}
\right.\right\}$$

 and 

$$ %\;\;\;\;\;\;\;\;\;\;\;\;\;\;\;\;\;\;\;\;\;\;\;\;\;\;\;\;\;\;\;\;\;
\#\left\{(\sigma_u,\tau)\left|
	      \begin{array}{l}
\bullet\ \sigma_u\in\Perm_d^{(\Delta_u)}\\
\bullet\ \tau\textnormal{ an (unmarked) transposition in }\Perm_d\\
\bullet\ \tau\circ\sigma_u=\tilde\sigma
             \end{array}
\right.\right\}.
$$
For the second factor, it is easier to multiply with $\tau=\tau^{-1}$ and count the number of transpositions $\tau$ satisfying $\tau\circ\tilde\sigma \in  \Perm_d^{(\Delta_u)}$. The requirement is satisfied if and only if both entries of $\tau$ come from the two different cycles of $\tilde\sigma$ which are joined to one cycle. We can thus choose one entry of the $m_1$ entries of one cycle, and one of the $m_2$ entries of the other, leading to $m_1\cdot m_2$ choices. Since this holds true for any $\tilde\sigma$, we can pull this factor in front of the sum. Our number then equals 
\begin{align*}(m_1&+m_2)\cdot \frac{1}{d!}\cdot \\ &\sum_{\tilde\sigma\in\Perm_d^{(\tilde\Delta)}} \#\left\{
(\sigma_v,\sigma_w)\left|
	    \begin{array}{l}
\bullet\ \sigma_v\in\Perm_d^{(\Delta_v)},\sigma_w\in\Perm_d^{(\Delta_w)}\\
\bullet\ \sigma_w\circ\sigma_v\circ\tilde\sigma=\operatorname{id}_{\Perm_d}\\
\bullet\ \langle\tilde\sigma,\sigma_v,\sigma_w\rangle_{\Perm_d}\textnormal{
acts transitively on }\{1,\ldots,d\}
             \end{array}
\right.\right\}.\end{align*}
The sum times $\frac{1}{d!}$ equals $H_d^g(\tilde\Delta,\Delta_v,\Delta_w)$ if $m_1\neq m_2$ and 
$\frac{1}{2}H_d^g(\tilde\Delta,\Delta_v,\Delta_w)$ if $m_1=m_2$ (because if $m_1=m_2$ there are two ways to mark the two preimages with ramification index $m_1=m_2$ above the point with ramification profile $\tilde\Delta$ which we count only once here since we have only a partially marked partition).
Since the dual graph has a wiener if and only if $m_1=m_2$ (leading to a factor of $\frac{1}{2}$ in the tropical multiplicity), the product equals the tropical multiplicity.

\vspace{0.4cm}

Now assume that the dual graph of the boundary point is as in case (2) of the proof of proposition
\ref{prop:1-1-correspondence}, i.e.\ as e.g.\ in the top row
on the left of Figure \ref{fig:tropResolutions}.
There is one component covering the other copy of $\PP^1$, and it
meets $C_1$ in one node. Then $\tilde \Delta$ is obtained from $\Delta_u$ by
summing two parts $m_1$ and $m_2$. 
Again, if we consider covers with simple ramification at $t$ close to $0$, we can count the ones which degenerate to this boundary point just as in equation \ref{eq-count}.

We claim that if
$\langle\tau,\sigma_u,\sigma_v,\sigma_w\rangle_{\Perm_d}$
acts transitively on $\{1,\ldots,d\}$ then so does
$\langle\tau\circ\sigma_u,\sigma_v,\sigma_w\rangle_{\Perm_d}$, hence we can drop the second transitivity condition.
Assume $\langle\tau,\sigma_u,\sigma_v,\sigma_w\rangle_{\Perm_d}$
acts transitively. For arbitrary $k,l\in\{1,\ldots,d\}$ we would like to have a word in
$\tau\circ\sigma_u$, $\sigma_v$, $\sigma_w$ and their inverses which as a permutation maps $k$
to $l$. Let $\tau$ be $(\tau_1,\tau_2)$. It joins two cycles $c_1$ and $c_2$ of
$\sigma_u$ (containing the elements $\tau_1$ and $\tau_2$ respectively) to a
cycle $c$ in $\tilde\sigma$ (obviously containing $\tau_1$ and $\tau_2$).
The remaining cycles are the same in both permutations. Therefore there are
$s,t\in\N$ such that $\tilde\sigma^s(\tau_1)=c^s(\tau_1)=\tau_2$ and
$\tilde\sigma^t(\tau_2)=c^t(\tau_2)=\tau_1$. 
Since $\langle\tau,\sigma_u,\sigma_v,\sigma_w\rangle_{\Perm_d}$ acts
transitively, we have a product $\delta_r\circ\ldots\circ\delta_1$ where each $\delta_i$ is one of the
permutations 
$\tau,\sigma_u,\sigma_v$ and $\sigma_w$ or their inverses, and which maps $k$ to $l$. Let
$k_i$ be $\delta_i\circ\ldots\circ\delta_1(k)$ for $i=1,\ldots,r$ and $k_0=k$.
Assume $\delta_i=\tau$ and $k_{i-1}$ is in the support of $\tau$. If
$k_{i-1}=\tau_1$ define $\delta_i'=\tilde\sigma^s$ and 
$\delta_i'=\tilde\sigma^t$ otherwise. Then clearly
$\delta_r\circ\ldots\circ\delta_i'\circ\ldots\circ\delta_1(k)=\delta_r\circ\ldots\circ\delta_i\circ\ldots\circ\delta_1(k)
$. Analogously
if $\delta_i=\sigma_u$ (or $\sigma_u^{-1}$) with $k_{i-1}$ in the support of $c_1$ or $c_2$, we can
substitute $\delta_i$ with powers of $\tilde\sigma$ (or $\tilde\sigma^{-1}$). In this way we produce the desired
word in the permutations in
$\langle\tau\circ\sigma_u,\sigma_v,\sigma_w\rangle_{\Perm_d}$ mapping $k$ to $l$.

After dropping the
second transitivity condition in equation \ref{eq-count}, we can as before write
the number as a product of two factors 
$$\frac{1}{d!}\cdot\sum_{\tilde\sigma\in\Perm_d^{(\tilde\Delta)}}\#\left\{
(\sigma_v , \sigma_w)\left|
	    \begin{array}{l}
\bullet\ \sigma_v\in\Perm_d^{(\Delta_v)},\sigma_w\in\Perm_d^{(\Delta_w)}\\
\bullet\
\tilde\sigma\circ\sigma_v\circ\sigma_w=\operatorname{id}_{\Perm_d}\\
\bullet\ \langle\tilde\sigma,\sigma_v,\sigma_w\rangle_{\Perm_d}\textnormal{
acts transitively on }\{1,\ldots,d\}
             \end{array}
\right.\right\}
$$
and
$$
\#\left\{(\sigma_u,\tau)\left|
	    \begin{array}{l}
\bullet\ \sigma_u\in\Perm_d^{(\Delta_u)}\\
\bullet\ \tau\textnormal{ an (unmarked) transpositions in }\Perm_d\\
\bullet\ \tau\circ\sigma_u=\tilde\sigma
             \end{array}
\right.\right\},$$
where we can pull the second factor out of the sum because it is the same for each $\tilde\sigma$: we count transpositions $\tau$ satisfying $\tau\circ\tilde\sigma\in \Perm_d^{(\Delta_u)}$. To obtain such a $\tau$, we can pick any entry in the joined cycle in $\tilde\sigma$, and pick as second entry one which is $m_1$ numbers away. If $m_1\neq m_2$, we have $m_1+m_2$ different choices. If $m_1=m_2$, we have $m_1$ choices but then for each choice two options for the marking of $\tau\circ\tilde\sigma$, so altogether we get $2m_1=m_1+m_2$ also. 

The first factor equals $H_d^g(\tilde\Delta,\Delta_v,\Delta_w)$. 
The product equals the tropical multiplicity.

\vspace{0.4cm}

Finally, assume that  the dual graph of the boundary point is as in case (3) of the proof of proposition
\ref{prop:1-1-correspondence}, i.e.\ as e.g.\ in the top row in the middle of Figure
\ref{fig:tropResolutions}.

There are two components $D_1$ and $D_2$ of genus $g_1$ and $g_2$ covering
the other copy of $\PP^1$ with degree $d_1$ and $d_2$ respectively, each meeting $C_1$ in one node. Then $\tilde\Delta$
is obtained from $\Delta_u$ by splitting the part $m$  into two parts $m_1$ and $m_2$. Moreover, $\tilde\Delta$ is naturally
divided into two partitions $\tilde\Delta^{(1)}$ and $\tilde\Delta^{(2)}$ of $d_1$ and $d_2$ respectively depending on whether the corresponding
ramification point is in $D_1$ or $D_2$. 
In the same way, the partitions
$\Delta_v$ and $\Delta_w$ are divided into
$\Delta_v^{(1)},\Delta_v^{(2)}$ and $\Delta_w^{(1)},\Delta_w^{(2)}$ respectively. 

If we consider covers with simple ramification at $t$ close to $0$, we can count the ones which degenerate to this boundary point as follows:
$$\frac{1}{d!}\cdot\sum_{S\subset\{1,\ldots,d\},\#S=d_1}
\#\left\{(\sigma_u,\sigma_v,\sigma_w,\tau)\right\},$$
where the tuples in the braces satisfy
\begin{enumerate}
 \item[(C1)] $\sigma_u\in\Perm_d^{(\Delta_u)}$,
\item[(d1)] $\sigma_v\in\Perm_d^{(\Delta_v) } ,\sigma_w\in\Perm_d^{(\Delta_w)}$;
\item[(C2)] $\tau$ an (unmarked) transposition in $\Perm_d$;
\item[(d2)] $ \tau\circ\sigma_u\circ\sigma_v\circ\sigma_w=\operatorname{id}_{\Perm_d}$;
\item[(d3)] $\langle\tau,\sigma_u,\sigma_v,\sigma_w\rangle_{\Perm_d}$ acts transitively on $\{1,\ldots,d\}$;
\item[(C3)] $\tau\circ\sigma_u=\tilde\sigma^{(1)}\circ\tilde\sigma^{(2)}$ (where
$\tilde\sigma^{(1)}$ and $\tilde\sigma^{(2)}$ are disjoint permutations acting
on the subset $S\subset \{1,\ldots,d\}$ resp.\ $S^c$ satisfying
$\tilde\sigma^{(1)}\in \Perm_S^{(\tilde\Delta^{(1)})}$ resp.\
$\tilde\sigma^{(2)}\in \Perm_{S^c}^{(\tilde\Delta^{(2)})}$);
\item[(A1)] there are permutations $\sigma_x^{(1)}\in\Perm_{S}^{(\Delta_x^{(1)})}$ for $x=v,w$ and
\item[(B1)] $\sigma_x^{(2)}\in\Perm_{S^c}^{(\Delta_x^{(2)})}$ for $x=v,w$ satisfying
\item[(d4)] $ \sigma_x^{(1)}\circ \sigma_x^{(2)}=\sigma_x$ for $x=v,w$;
\item[(A2)] $\tilde\sigma^{(1)}\circ\sigma_v^{(1)}\circ\sigma_w^{(1)}=\operatorname{id}_{\Perm_S}$;
\item[(B2)] $\tilde\sigma^{(2)}\circ\sigma_v^{(2)}\circ\sigma_w^{(2)}=\operatorname{id}_{\Perm_{S^c}}$;
\item[(A3)] $\langle\tilde\sigma^{(1)},\sigma_v^{(1)},\sigma_w^{(1)}\rangle_{\Perm_{S}}$ acts transitively on $S$;
\item[(B3)] $\langle\tilde\sigma^{(2)},\sigma_v^{(2)},\sigma_w^{(2)}\rangle_{\Perm_{S^c}}$ acts transitively on $S^c$.
\end{enumerate}
Due to (C3) and (C1), $\tau$ must have one entry in $S$ and one in $S^c$, so the
transitivity  condition (d3) is implied by (A3) and (B3). Moreover, (B2) and
(A2) imply (d2).
Below, we count the possibilities for tuples  $(\sigma_v^{(i)},\sigma_w^{(i)})$ for $i=1,2$. The permutation $\sigma_x$ is then by (d4) given as the product of the two entries and we can thus neglect it and condition (d1) which is implied by (A1) and (B1).
Finally, ordering the tuples by the different possibilities
for $\tilde\sigma^{(1)}$ and $\tilde\sigma^{(2)}$ we can write the above number as

\begin{align*}\frac{1}{d!}\cdot\sum_{S\subset\{1,\ldots,d\},\#S=d_1}\sum_{
\tilde\sigma^{(1)} \in\Perm_{S}^{
(\tilde\Delta^{(1)})}}&\sum_{\tilde\sigma^{(2)}\in\Perm_{S^c}^{(\tilde\Delta^{(2)})}}
\#\left\{(\sigma_v^{(1)},\sigma_w^{(1)})\left|\right.\mbox{(A)} 
\right\}\cdot \\&
\#\left\{(\sigma_v^{(2)},\sigma_w^{(2)})\left|\right.\mbox{(B)} 
\right\}\cdot 
\#\left\{(\sigma_u,\tau)\left|\right.\mbox{(C)} 
\right\},\end{align*}

where a capital letter stands for the three conditions labeled accordingly.

 The last
factor in each summand equals $m_1\cdot m_2$ for all choices of the $\tilde\sigma^{(i)}$ by the same argument as in the first case.
Instead of summing over all subsets of $\Perm_d$ of size $d_1$ we can fix
without restriction $S=\{1,\ldots,d_1\}$ and multiply by $\binom{d}{d_1}$. 
Furthermore, the two factors in each summand each depend on only one summation index, so we can sort the sums accordingly. Notice that 
$$\sum_{\tilde\sigma^{(1)}\in\Perm_{S}^{(\tilde\Delta^{(1)})}}
\#\left\{(\sigma_v^{(1)},\sigma_w^{(1)})\big|\mbox{(A)}
\right\} = d_1!\cdot H_{d_1}^{g_1}
(\tilde\Delta^{(1)},\Delta_v^{(1)},\Delta_w^{(1)})\
$$
and
 $$\sum_{\tilde\sigma^{(2)}\in\Perm_{S^c}^{(\tilde\Delta^{(2)})}}
 \#\left\{(\sigma_v^{(2)},\sigma_w^{(2)})\left|\mbox{(B)}
 \right.\right\}) = d_2!\cdot H_{d_2}^{g_2}
 (\tilde\Delta^{(2)},\Delta_v^{(2)},\Delta_w^{(2)}),
 $$

so we get
\begin{eqnarray*}
 &&m_1\cdot
m_2\cdot\frac{1}{d!}\cdot\binom{d}{d_1}d_1!\cdot H_{d_1}^{g_1}
(\tilde\Delta^{(1)},\Delta_v^{(1)},\Delta_w^{(1)})\cdot d_2!\cdot H_{d_2}^{g_2}
(\tilde\Delta^{(2)},\Delta_v^{(2)},\Delta_w^{(2)})\\
&=&m_1\cdot
m_2\cdot
H_{d_1}^{g_1}(\tilde\Delta^{(1)},\Delta_v^{(1)},\Delta_w^{(1)})\cdot H_{d_2}^{
g_2}
(\tilde\Delta^{(2)},\Delta_v^{(2)},\Delta_w^{(2)}).
\end{eqnarray*}
This equals the tropical multiplicity.

\end{proof}

\begin{cor}\label{cor-one}
 The degree of the tropical branch map $\br^{\trop}:\Mg(\CL,\Delta)\rightarrow \CL $ from a \emph{one}-dimensional space $\Mg(\CL,\Delta)$ (i.e.\ $2g-2-d+\#\Delta=1$) is constant.
In particular, if we consider all resolutions of the star-shaped cover and group their multiplicities into three sums corresponding to the three ends of $\CL$ to which the trivalent vertex can be mapped, the three sums agree.
\end{cor}
\begin{proof}
 This follows from lemma \ref{lem:1:1corrWithMult} and proposition
\ref{prop:brIsCover}.
\end{proof}
\begin{example}
 If we add the tropical multiplicities for each column in figure
\ref{fig:tropResolutions}, we get 
\begin{align*}&2\cdot 4 \cdot H^1_5((4,1),(5),(3,2)) + 2 \cdot H^1_5((2,3),(5),(2,3))\\&+ 2\cdot H^0_5((1,1,1,2),(5),(2,3))\\=&2\cdot4\cdot2+2\cdot1+2\cdot6=30\end{align*} for the
left column (the first factor $2$ comes from the fact that the upper left
figure stands for two different types of cover due to the different possibilities
to mark the weight-$1$-edges over $u$).
This is the sum of tropical multiplicities of resolutions where the trivalent vertex is mapped to $u$.
In the same way, we get
\begin{align*}&2\cdot 3\cdot H^0_2((1,1),(2),(2))\cdot H^1_3((3),(3),(3))\\&+
2\cdot 3\cdot H^0_5((3,1,1),(2,3),(2,3))+ 4\cdot H^0_5((3,1,1),(4,1),(2,3))\\=&
(2\cdot3)\cdot(1\cdot\frac{1}{3})+(2\cdot3)\cdot2+4\cdot4=30\end{align*} for the middle
column corresponding to resolutions where the trivalent vertex is mapped to $v$
and \begin{align*}&5\cdot H^1_5((3,1,1),(5),(5))+ \frac{1}{2}\cdot H^0_5((3,1,1),(5),(3,1,1))\\&+ 2\cdot H^0_5((3,1,1),(5),(2,2,1))\\=&
5\cdot4+\frac{1}{2}\cdot4+2\cdot4=30\end{align*} for the right column corresponding to resolutions where the trivalent vertex is mapped to $w$.
\end{example}

\subsection{The proof of Theorem \ref{thm}}
Using the duality in the one-dimensional case, we can now prove Theorem
\ref{thm}. So assume we are given $\Delta,d$ and $g$ such that
$2g-2-d+\#\Delta\geq0$.

\begin{proof}[Proof of Theorem \ref{thm}]
We refine the fan $\CL^r$ by adding the diagonals
$D_{ij}$ defined as $\{(p_1,\ldots,p_r)\ |\ p_k\neq c\ \forall k=1,\ldots,r,\
p_i=p_j\}$ for $i\neq j$ as codimension-1-faces, where $c$ as before denotes the center of the line $\CL$. Let us call the new fan by
abuse of notation $\CL^r$ as well. The point configurations in the interior of top-dimensional faces of $\CL^r$ are in general position. The degree of $\br^{\trop}$
is constant on any top-dimensional face, since the preimages of two different
point configurations in the same face contain the same combinatorial types.

As $\CL^r$ is connected in codimension $1$ it is sufficient to see that the
degree of $\br^{\trop}$ does not change if we cross a
codimension-1-face in $\CL^r$.

Let us first assume that we cross a diagonal, that is beginning from a point configuration $P$
in general position two branch points on one of the ends of $\CL$ change their positions. We call the new point configuration $P'$. One can see easily that we have exactly the same combinatorial types of curves in the preimages of $P$ and $P'$, they just differ by their vertex labelings (see also  \cite{CJM10}, Lemma 5.27). 
Thus the degree of $\br^{\trop}$ is constant when crossing this diagonal.

%  The combinatorial type $\alpha$ of
% any fixed curve in the preimage of $P$ determines a partition $\Delta_\alpha$ whose parts are
% the weights of the edges $e_1,\ldots,e_t$ mapping to the $u$-end and adjacent to a vertex over the center.
% For a point $P'$ in $\CL^r$ close to $P$, the preimages under $\br$ are resolutions of the preimages of $P$. We can interpret each 
% 
% 
% The contribution of $\alpha$ to the degree of $\br$
% is the product of the contri the double Hurwitz number $H^d_{\tilde g}(\Delta_u,\Delta_\alpha)$, where $\tilde g$ is the genus of the part mapping to $u$ of the underlying graph after cutting the edges $e_1,\ldots,e_t$ and factors that come from the remaining part of the cover. The Hurwitz number does not
% depend on the position of the branch points (see \cite{CJM10}). 

% , the degree of
% $\br$ does not change when crossing a diagonal.

Now let us fix a point configuration $P$ on a codimension-1-face in $\CL^r$ which is not a diagonal, that is a point configuration
where exactly one point is the center $c$ of $\CL$. 
The combinatorial types of the preimages with
respect to $\br^{\trop}$ have exactly one simple ramification over the center and all
other simple ramifications over the ends. For a fixed type $\alpha$ the
top-dimensional cones adjacent to $D_\alpha$ in $\Mg(\CL,\Delta)$ correspond to the resolutions of the simple ramification over the center
as described in section \ref{one-dimCase}. We can thus interpret their contribution to the degree of $\br^{\trop}$ as a product of a local factor corresponding to the one-dimensional resolution and factors from the remaining parts of the cover, which are the same in any case.
Since by corollary \ref{cor-one} the local factors add to a contribution which does not depend on the end of $\CL$ above which we pull the simple ramification, the degree of $\br^{\trop}$ is constant locally around $P$.

% The contribution of these resolutions
% with an additional branch point on the $u$-end to the degree of
% $\br$ consists of the corresponding tropical triple Hurwitz number and factors
% from the remaining part of the curve. The latter is the same for all
% resolutions and the first factor does not not depend on whether we resolve in
% $u$-,$v$- or $w$-direction by lemma \ref{lem:1:1corrWithMult} and proposition
% \ref{prop:brIsCover}.
\end{proof}

\bibliographystyle{plain}
\bibliography{bibliographie}

\end {document}